\documentclass{elsart3-1}

\usepackage{amsmath}
\usepackage{amssymb}
\usepackage[english,francais]{babel}

\def\KK{{\mathbb K}}  \def\AA{{\mathbb A}}

  \def\PP{{\mathbb P}}
 \def\QQ{{\mathbb Q}} 
\def\ZZ{{\mathbb Z}}

\def\Res{{\mathrm{Res}}}
\def\a{{\mathbf{a}}}

\begin{document}

\begin{frontmatter}

\selectlanguage{english}
\title{On the irreducibility of multivariate subresultants}

\vspace{-2.6cm}

\selectlanguage{francais}
\title{Sur l'irr\'eductibilit\'e des sous-r\'esultants multivari\'es}

\selectlanguage{english}
\author[authorlabel1]{Laurent Bus\'e}
\ead{lbuse@sophia.inria.fr}
\author[authorlabel2]{Carlos D' Andrea}
\ead{cdandrea@math.berkeley.edu}

\address[authorlabel1]{INRIA, GALAAD, 2004 route des Lucioles, B.P. 93,
  06902 Sophia-Antipolis cedex, France. T\'el:
  +33-(0)4-92-38-76-51, fax: +33-(0)4-92-38-79-78.}
\address[authorlabel2]{Department of Mathematics UC Berkeley
 \&   The Miller Institute for
Basic Research in Science. 1089 Evans Hall, Berkeley, CA 94720-3840
USA. T\'el: +1-(510)-642-6947, fax: +1-(510)-642-8204.}

\begin{abstract}
 Let $P_1,\ldots,P_n$ be generic homogeneous polynomials in $n$ variables of degrees $d_1,\ldots,d_n$ respectively.
We prove that if $\nu$ is an integer satisfying
${\sum_{i=1}^n d_i}-n+1-\min\{d_i\}<\nu,$
then all multivariate subresultants associated to the family $P_1,\ldots,P_n$ in degree $\nu$ are irreducible. We show that the lower bound is sharp.
As a byproduct, we get a formula for computing the residual resultant of $\binom{\rho-\nu +n-1}{n-1}$ smooth isolated points in $\PP^{n-1}.$

\vskip 0.5
\baselineskip

\selectlanguage{francais}
\noindent{\bf R\'esum\'e}
\vskip 0.5\baselineskip
\noindent
Soient $P_1,\ldots,P_n$ des polyn\^omes homog\`enes g\'en\'eriques en
$n$ variables de degr\'e respectif $d_1,\ldots,d_n$. Nous montrons que
si $\nu$ est un entier tel que ${\sum_{i=1}^n
  d_i}-n+1-\min\{d_i\}<\nu$,  tous les sous-r\'esultants multivari\'es
de degr\'e $\nu$ des polyn\^omes $P_1,\ldots,P_n$ sont
 irr\'eductibles. Nous montrons \'egalement que cette borne est
 atteinte dans des cas particuliers. Comme cons\'equence directe nous obtenons
 une nouvelle formule pour le calcul du r\'esultant
 r\'esiduel de $\binom{\rho-\nu +n-1}{n-1}$ points lisses isol\'es
 dans $\PP^{n-1}$.

\end{abstract}
\end{frontmatter}

\maketitle

\selectlanguage{english}
Classical subresultants of two univariate polynomials have been studied by Sylvester in the foundational work
\cite{syl}.
Multivariate subresultants, introduced in \cite{cha1}, provide a criterion for over-constrained
polynomial systems to have Hilbert function of prescribed value,
generalizing the classical case.
To be more precise, let $\KK$ be a field. If $P_1,\dots,P_s$ are homogeneous polynomials in $\KK[X_1,\dots,X_n]$ with $d_i=\deg(P_i)$ and $s\leq n,$
$H_{d_1,\dots,d_s}({\bf .})$ is the Hilbert function of a complete intersection given by $s$ homogeneous polynomials in $n$ variables of degrees $d_1,\dots,d_s,$
and $S$ is a set of $H_{d_1,\dots,d_s}(\nu)$ monomials of degree $\nu,$ the subresultant $\Delta_S^\nu$ is a polynomial in the coefficients of the $P_i's$ of degree
$H_{d_1,\dots,d_{i-1},d_{i+1},\dots,d_n}(\nu-d_i)$ in the coefficients of $P_i$ ($i=1,\dots,s)$ having
the following universal property: $\Delta_S^\nu\neq0$ if and only if $I_\nu+\KK\langle S\rangle = \KK[X_1,\dots,X_n]_\nu,$
where $I_\nu$ is the degree $\nu$ part of the ideal generated by the
$P_i$'s (see \cite{cha1}).

Multivariate subresultants have been used in computational algebra for polynomial system solving (\cite{GLV3},\cite{sza2}) as well as for providing
explicit formulas for the representation of rational functions
(\cite{jou3,DJ,DK,mul}). The study of their properties is an
active area of research (\cite{cha2,cha3,DJ,DK,elk}). In particular, it is important to know which $S$ verify $\Delta_S^\nu\neq0$, and which of these $\Delta_S^\nu$ are irreducible (see the final remarks and
open questions in \cite{cha1} and the conjectures in \cite{DK}). Partial results have been obtained in this direction. In \cite{cha4} it is shown that, if $s=n$ and
$
{\sum_{i=1}^n d_i} -n-\min\{d_i\}< \nu,
$
then for every set $S$ of monomials
of degree $\nu$ and cardinal $H_{d_1,\dots,d_n}(\nu),$ the polynomial $\Delta_S^\nu$ is not identically zero. Moreover, in \cite{cha3},
it is also proved that if
$s=n,\,\nu=\sum_{i=1}^n d_i-n,$ and $S=\{x_j^\nu\}$ for $j=1,\ldots,n,$ then $\Delta_S^\nu$ is an irreducible polynomial in the coefficients of the $P_i's.$
In \cite[Lemma 4.2]{elk} the irreducibility of $\Delta_S^\nu$ is shown for $s=n=2,\,\max\{d_1,d_2\}\leq\nu,$ and $S=\{X_2^\nu,X_1X_2^{\nu-1},\dots,
X_1^{H_{d_1,d_2}(\nu)-1}X_2^{\nu-H_{d_1,d_2}(\nu)+1}\}.$

In this note we study the irreducibility problem in
the case $s=n$. Let us introduce some notations in order to state our
result.
Let $\rho:=\sum_{i=1}^n (d_i-1)$. For $i=1,\ldots,n$ and $\alpha\in \ZZ_{\geq0}^n$ such that $|\alpha|=d_i,$ introduce a new variable $c_{i,\alpha}.$
Let $\AA:=\ZZ\left[c_{i,\alpha},i=1,\ldots,n,\,|\alpha|=d_i\right]$
and
set
\begin{equation}\label{input}
P_i(x_1,\ldots,x_n):=\sum_{|\alpha|=d_i}c_{i,\alpha}x^\alpha.
\end{equation}
\noindent {\bf Theorem} {\it
For every $\nu$ such that
$\rho-\min\{d_i\}+1<\nu$ and every set $S$ of monomials of degree
$\nu$ and cardinality $H_{d_1,\dots,d_n}(\nu),$ the subresultant $\Delta_S^\nu(P_1,\ldots,P_n)$ is irreducible in $\AA.$}
\smallskip

Observe that, if $n=2,$ then $\rho-\min\{d_i\}+1=d_1+d_2-2-\min\{d_i\}+1=\max\{d_i\}-1,$ and this is equivalent to $\max\{d_i\}\leq\nu,$
so our result contains those in \cite{elk}.

\medskip

\noindent {\bf Proof of the Theorem:}
For simplicity we assume hereafter that $d_1\geq \ldots \geq d_n\geq 1$. First observe that if $\nu > \rho$ then $\Delta_S^\nu$ is simply a resultant,
and is hence known to be irreducible. So, we can suppose w.l.o.g. that $d_n>1.$ We thus only have to consider integers $\nu$ such that
\begin{equation}\label{nu}
  \rho \geq \nu > \rho - d_n+1=\sum_{i=1}^{n-1}(d_i-1),
  \end{equation}
where we recall that $\rho=\sum_{i=1}^n(d_i-1)$. We begin by computing
the multi-degree of the subresultants $\Delta_S^\nu$; we know (see
\cite{cha1}) that
$$\deg_{P_i}(\Delta_S^\nu)=H_{d_1,\ldots,d_{i-1},d_{i+1},\ldots,d_n}
(\nu-d_i).$$
But from the standard short exact sequence
{\scriptsize
$$0 \rightarrow \frac{R}{(f_1,\ldots,f_{i-1},f_{i+1},\ldots,f_{n})}(-d_i) \xrightarrow{\times f_i}
\frac{R}{(f_1,\ldots,f_{i-1},f_{i+1},\ldots,f_{n-1})} \rightarrow
\frac{R}{(f_1,\ldots,f_{n})} \rightarrow 0,$$}

\noindent where $f_1,\ldots,f_n$  are homogeneous polynomials of respective
degree $d_i$ in a graded polynomial ring $R$ and $f_1,\ldots,f_n$ is a
complete intersection in $R$, we deduce
$$H_{d_1,\ldots,d_{i-1},d_{i+1},\ldots,d_n}(t-d_i)=H_{d_1,\ldots,d_{i-1},d_{i+1},\ldots,d_n}(t)-H_{d_1,\ldots,d_n}(t)$$
for all integer $t$. It follows that for all integer $\nu\geq \rho-d_n+1$,
\begin{equation}\label{degsubres}
\deg_{P_i}(\Delta_S^\nu)=\frac{d_1\ldots
  d_n}{d_i}-H_{d_1,\ldots,d_n}(\nu)=\frac{d_1\ldots
  d_n}{d_i}-\binom{\rho-\nu+n-1}{n-1},
\end{equation}
where that last equality comes from the facts that
$H_{d_1,\ldots,d_n}(\rho-t)=H_{d_1,\ldots,d_n}(t)$ for all integer
$t$, and $H_{d_1,\ldots,d_n}(t)=\binom{t+n-1}{n-1}$ for all $0
\leq t<d_n$. We define ${\a}:=\binom{\rho-\nu+n-1}{n-1}.$ As
$\a$ does not depend on $i \in \{1,\ldots,n\}$ and residual (or reduced) resultants of $\a$ isolated points in $\PP^{n-1}$ have the same degree in the
coefficients of $P_i$ as the right hand side of (\ref{degsubres}), this suggest that we compare $\Delta_S^\nu$ with residual resultants.

We will work with  an ideal $G$ defining $\a$ points in $\PP^{n-1}$
which is generated in degree at most $d_n$ and such that $G_{d_n-1}\neq 0$.
Ideals defining $\a$ points in sufficiently generic position are generated in degree exactly $\rho-\nu+1$
(see \cite[Proposition 4]{GO}). Since by \eqref{nu} we have $d_n > \rho -\nu +1$, we thus choose such an ideal $G=(g_1,\ldots,g_m)$,
where $\deg(g_i)=\rho-\nu+1$ for all $i=1,\ldots,m$, defining $\a$ points in generic position (see \cite{GO} for the definition of ``generic position''), and hence locally a
complete intersection.

Now consider the following specialization of polynomials $P_i$'s
\begin{equation}\label{tg}
P_i\mapsto \overline{P}_i:=\sum_{j=1}^{m} p_{ij}(x)g_j(x),
\end{equation}
where $p_{ij}(x)=\sum_{|\alpha|=d_i-\rho+\nu-1}c_{ij}^{|\alpha|}x^\alpha$ is a
generic polynomial of degree $d_i-\rho+\nu-1$. There exists a resultant associated to the system $\overline{P}_1,\ldots,\overline{P}_n$, called the \emph{residual resultant}. We denote it by $\Res_G(\overline{P}_1,\ldots,\overline{P}_n)$. Let us recall its main properties (see  \cite{bus} \S3.1).

\begin{itemize}
\item $\Res_G(\overline{P}_1,\ldots,\overline{P}_n)$ is a homogeneous and \emph{irreducible} polynomial in the ring of all the coefficients $\QQ[c_{ij}^{|\alpha|}]$,\\
\item For any given specialization of the coefficients $c_{ij}^{|\alpha|}$'s sending $\overline{P}_i$ to $Q_i$, we have

$$\Res_{G}(Q_{1},\ldots,Q_{n})=0  \ \mbox{if and only if} \   
 (Q_1,\ldots,Q_n)^{sat} \varsubsetneq G=G^{sat},$$
\item $\Res_G(\overline{P}_1,\ldots,\overline{P}_n)$ is multi-homogeneous:
it is homogeneous in the coefficients of each polynomials $\overline{P}_i$, $i=1,\ldots,n$, and we have
$$\deg_{\overline{P}_i}(\Res_G(\overline{P}_1,\ldots,\overline{P}_n))=\frac{d_1\ldots
  d_n}{d_i}-\a.$$
\end{itemize}

We are now going to  compare this residual resultant with the
specialized subresultant
$\Delta_S^\nu(\overline{P}_1,\ldots,\overline{P}_n)$, which is non-zero as proved in \cite{cha3}. We claim that we have the following implications:

\begin{equation}\label{impl}
\Delta_S^\nu(Q_1,\ldots,Q_n)\neq 0 \Rightarrow
H_{(\underline{Q})}(\nu)=\a \Rightarrow H_{(\underline{Q})}(t)=\a \
\mbox{for all} \ t \geq \nu
\Rightarrow  \Res_G(Q_1,\ldots,Q_n)\neq 0,
\end{equation}

\noindent where $H_{(\underline{Q})}(.)$ denotes the Hilbert function associated to the ideal $(Q_1,\ldots,Q_n)$. Only the second implication needs to be proved, the others follow directly from
the algebraic properties of resultants and subresultants.
We know that $H_G(t)=\a$ for all $t\geq \rho-\nu+1$ (see \cite{GO}), and since we have supposed \eqref{nu}, it is a straightforward computation to show that
 $\nu\geq \rho-\nu+1$. It follows that, by hypothesis, the ideals $G$ and
 $(\underline{Q})$ coincide in degree $\nu$ and 
 have Hilbert function $\a$ in this degree. As they are both generated 
 in degree at most $\nu$ this implies that they coincide in all higher 
 degrees, and therefore they both have Hilbert function equal to $\a$
 in these degrees, because $G$ is the defining ideal of a set of points.

Due to (\ref{impl}) and the irreducibility of the residual resultant, we deduce that  $\Res_G(\overline{P}_1,\ldots,\overline{P}_n)$ divides $\Delta_S^\nu(\overline{P}_1,\ldots,\overline{P}_n)$.
But both polynomials have the same degree, so they
must be equal up to a rational number (giving a new formula for
computing this residual resultant using \cite{cha2}). Since this residual resultant
is irreducible, and since $\Delta_S^\nu$ and $\Delta_S^\nu(\overline{P}_1,\ldots,\overline{P}_n)$
have the same multi-degree, this shows that $\Delta_S^\nu$ is irreducible in $\QQ[\mathrm{coeff}(P_i)]$.

It remains to prove that $\Delta_S^\nu$ is irreducible in
$\ZZ[\mathrm{coeff}(P_i)].$ As it is irreducible in $\QQ[\mathrm{coeff}(P_i)],$ we only have to show that $\Delta_S^\nu$ has content $\pm1.$ Suppose that
this is not the case, and let $p\in\ZZ$ be a prime dividing the content of $\Delta_S^\nu.$ Let $k$ be the algebraic closure of $\ZZ_p.$ This implies that
$\Delta_S^\nu=0$ in $K:=k(\mathrm{coeff}(P_i)),$ and hence $S$ is linearly dependent in $K[x_1,\dots,x_n] / \langle P_1,\dots,P_n \rangle,$ contradicting the
main result of \cite{cha3}.

\bigskip
\noindent {\bf Reducibility in lower degrees:}
We now exhibit some sets $S$ of degree  $\nu=\rho-\min\{d_i\}+1$ such
that $\Delta_S^\nu$ factorizes. This shows that the lower bound in our
theorem  is sharp.

\smallskip
\begin{itemize}
\item[$\bullet$] {\bf $\mathbf{n=2,\,d_1>d_2}$:} In this case,
  $\nu=d_1-1\geq d_2,$ and $H_{d_1,d_2}(\nu)=d_2.$ Thus $\Delta_S^\nu$ 
  can be here computed with Sylvester type matrices
  \cite{syl}. However, setting
  $f_2=c_0x_1^{d_2}+c_1x_1^{d_2-1}x_2+\cdots+c_{d_2}x_2^{d_2}$, the universal property of the subresultant
  $\Delta_S^\nu$ shows immediatly that it is a power of $c_0$, and we
  have already seen that its degree is $d_1-d_2+1$; it follows that
  $\Delta_S^\nu=c_0^{d_1-d_2+1}$, so it can not be irreducible.

  \smallskip
\item[$\bullet$]{\bf $\mathbf{n>2,\ d_1-1>d_2=d_3=\ldots=d_n=1}$:} Again in this case, $\nu=d_1-1$ and $H_{d_1,d_2}(\nu)=1.$ Choose $S=\{x_1^{\nu}\}$ and, if
$f_i=c_{1i}x_1+\ldots+c_{ni}x_n,\,i=2,\ldots,n,$ we set $\delta:=\det\left(c_{ij}\right)_{2\leq i,j\leq n}.$
Applying Lemma $4.4$ in \cite{DJ} to this situation, we get that
$\Delta_S^\nu=\delta^{\nu}.$ So, $\Delta_S^\nu$ is not irreducible.
\end{itemize}

\smallskip
\noindent {\bf Acknowledgements:} This work started during
the special year (2002--2003) on ``Commutative Algebra'' at the Mathematical Sciences Research Institute (MSRI Berkeley).
We wish to thank MSRI for its hospitality and for the wonderful
working atmosphere it provided. We are also grateful to David Cox,
Gabriela Jeronimo and Bernd Sturmfels for their helpful comments on
preliminary versions of this paper, and to the anonymous referee for 
useful remarks improving the presentation of this work.

\end{document}